\documentclass[12pt]{article}
\title{On simultaneous arithmetic progressions on elliptic curves}
\author{Irene Garc\'{\i}a--Selfa \and Jos\'e M. Tornero}
\date{February, 2006}

\newcommand{\QQ}{{\mathbf Q}}
\newcommand{\vs}{\vspace{.15cm}}

\begin{document}

\maketitle

\abstract{In this paper we study elliptic curves which have a
number of points whose coordinates are in arithmetic progression.
We first motivate this diophantine problem, prove some results, 
provide a number of interesting examples and, finally point out open
questions which focus on the most interesting aspects of the problem for us.}

\vs\vs


\vs


\section{Introduction}

During this paper, all standard results
unless otherwise explicitely stated are taken from \cite{Silv}.
Extensive use has been made of Maple V + APECS (by Ian Connell, 
McGill University).

We will deal with elliptic curves defined over a field $K$ by a
Weierstrass equation, that is
$$
E: Y^2 + a_1XY + a_3Y = X^3 + a_2X^2 + a_4X + a_6, \quad a_i \in K.
$$

We will denote, as usual, $E(K)$ the locus of the above
equation, together with the point at infinity, $O=(0:1:0)$.

Changes of variables preserving this form are those given by
$$
X' = u^2X + r, \quad Y' = u^3Y + sX + t;
$$
\noindent and we will consider that two equations related by such
a change of variables represent the same curve (equivalently, we
will deal with elliptic curves up to so--called Weierstrass changes of
variables).

Consider then $P_0,...,P_n \in E(K)$, with $P_i = (x_i,y_i)$
such that $x_0,...,x_n$ is an arithmetic progression. Then we say
that $P_0,...,P_n$ are in $x$--arithmetic progression ($x$--a.p.)
and also $E$ is said  to have an $x$--arithmetic progression of
length $n+1$. From the previous remarks, this does not depend on
the Weierstrass equation considered.

The same definition goes for $y$--arithmetic progressions
($y$--a.p.). However, in this case, changes of variables (even
those which preserve Weierstrass equations) can create and destroy
$y$--arithmetic progressions.

\vs\vs

\noindent
{\bf Example.--} Let us consider the equation over $\QQ$
$$
E: Y^2 - \frac{5}{16}XY + \frac{1}{64} Y = X^3 - \frac{1}{64} X^2,
$$
which verifies that
$$
\left( \frac{1}{8}, \frac{-4}{128} \right), \;
\left( \frac{-1}{32}, \frac{-3}{128} \right), \;
\left( \frac{5}{64}, \frac{-2}{128} \right), \;
\left( \frac{1}{32}, \frac{-1}{128} \right), \;
$$
$$
\left( \frac{1}{64}, 0 \right), \;
\left( \frac{3}{64}, \frac{1}{128} \right), \;
\left( \frac{1}{16}, \frac{2}{128} \right)
\in E(\QQ).
$$

The reader can easily check that, after the change of variables
$$
Y' = Y + X,
$$
\noindent for instance, the corresponding points are not in
$y$--a.p. Hence we can properly talk of $x$-a.p. in a curve, but
if we speak of $y$-a.p. in a curve we must bear in mind that we are
considering a specific equation.  

\vs\vs

This paper studies with elliptic curves which have a simultaneous
arithmetic progression. First we need a proper definition of these
progressions. Let us consider $P_0,...,P_n$ as above. If we ask
for both $x_0,...,x_n$ and $y_0,...,y_n$ to be arithmetic
progressions then the problem is far too easy, as $P_0,...,P_n$
must be collinear and hence $n \leq 2$. Lots of examples can be
found with this property; for instance all curves in the family
$$
E(b): Y^2 + (2b-1)XY + bY = X^3 - bX^2
$$
\noindent have the arithmetic progression $(0,-b),\; (b,0),\; (2b,b)$.

\vs\vs

\noindent {\bf Definition.--} With the above notation,
$P_0,...,P_n$ are a simultaneous arithmetic progression (s.a.p.)
or the equation $E$ is said to have a simultaneous arithmetic
progression if:

\begin{enumerate}
\item [(a)] $x_0,...,x_n$ are in arithmetic progression (called the support
of the s.a.p.).

\item [(b)] There exists a permutation $\sigma$ in the symmetric group of 
$n+1$ elements $S_{n+1}$ such that $y_{\sigma(0)},...,y_{\sigma(n)}$ are 
in arithmetic progression.
\end{enumerate}

The definition is clearly symmetric: it is equivalent (up to point
ordering) to saying that $y_0,...,y_n$ are in arithmetic
progression and there exists $\mu \in S_{n+1}$ such that
$x_{\mu(0)},...,x_{\mu(n)}$ are in arithmetic progression, but
this version showed more useful for computational purposes.

\vs\vs

With this definition, at least three straight problems arise:

\begin{enumerate}
\item [(a)] The detection problem: Given an elliptic curve, does
there exist an algorithm for deciding whether it contains or not a s.a.p. of length $n$
(giving as well a change of variables if needed)?

\item [(b)] The subsequence problem: If an elliptic curve has a
s.a.p. of length $n+1$, does it possess a s.a.p. of length $n$?
(Note this is not at all trivial from the definition). 

\item [(c)] The bound problem: Is there a bound for the possible lengths of 
s.a.p. in elliptic curves?
\end{enumerate}

Trying to understand these three problems, we have developed some
computational methods (actually, two partial answer to the detection problem) 
whose application may shed some light.
Nevertheless the results achieved can be considered only as a
first step towards a fully satisfactory understanding of these
sequences. We have managed to prove the following results:

\vs\vs

\noindent {\bf Theorem 1.--} Given an elliptic curve with an
$x$--a.p., there exists an algorithm which decides whether or not the curve
also has a s.a.p. with the given $x$--a.p. as support.

\vs\vs

\noindent {\bf Theorem 2.--} There are integers $n$ such that there are examples of 
s.a.p. of length $n$ which do not contain any s.a.p. of length $n-1$.

\vs\vs

\noindent {\bf Theorem 3.--} There are no elliptic curves defined over $\QQ$ 
with s.a.p. of length $7$. There are only finitely many non--isomorphic 
curves defined over $\QQ$ with s.a.p. of length $6$.

\vs\vs

We will finish this introduction with a word on motivation. At
first our interest was drawn to this subject by the articles of
Bremner--Silverman--Tzanakis (\cite{BST}) and Bremner
(\cite{BEM}). Apparently these papers had their starting point in
the featuring of $x$--a.p. as by--product of a latin square
problem (see more on this in \cite{BAA1,BAA2}). However, highly
interesting results were sketched in both papers around the
relationship between the existence of arithmetic progressions on a
certain elliptic curve and its rank. In this same line a paper by
Campbell appeared (\cite{Campbell}) pointing out far--reaching
questions, probably too difficult for the state--of--the--art.

The history of the problem, though, can be traced back quite
further, as (for the specific case of Mordell curves) it was
treated previously by S.P. Mohanty (\cite{Mohanty}) who studied
$x$ and $y$--a.p. separatedly and by Lee and V\'elez (\cite{LV})
who first treated s.a.p., if only in the naive form mentioned
above, without permutations involved. The motivation for these 
first works was, as many other times in the history of number 
theory, purely diophantine.

We became interested in this specific problem while trying to improve Bremner's
record of longest $x$--a.p. by narrowing the search\footnote{The longest 
$x$--a.p. found in an elliptic curve has $8$ terms (\cite{BEM}); for $y$--a.p. 
the record is $7$ so far (\cite{GT}).}. Our first attempts were shown in \cite{GT}, 
using a specific kind of s.a.p., which allowed us to find examples of s.a.p. of 
length $5$.  These methods were not at all exhaustive, as it accurately pointed 
out by Bremner in his MathSciNet review. After this work, we feel that some of 
the problems posed are worth a closer look and the setup remains challenging. 
As Bremner points out in \cite{BEM}: {\em ``Questions in number theory that 
interrelate two group structures are easily posed, but often lead to intractable problems''}.

\section{The detection problem}

Let us consider a set of points $P_0=(x_0, y_0),...,P_n=(x_n,
y_n)$ in an elliptic curve, defined over $K$ by a Weierstrass
equation:
$$
E:Y^2 + a_1XY + a_3Y = X^3 + a_2X^2 + a_4 X + a_6, \ \ \mbox{with } a_i \in K.
$$

Let us suppose the points $P_0,...,P_n \in E(K)$ to be in $x$--a.p.
We are interested then on finding, if there exists any, a change
of variables, preserving the Weierstrass form of $E$, which
transforms $P_i=(x_i, y_i)$ into $P'_i=(x'_i, y'_i)$ so that
$P'_0,...,P'_n$ is a s.a.p. on the corresponding equation $E'$.
This change of variables must be of the form
$$
X' = u^2X+s, \ \ \  Y' = u^3Y+rX+t.
$$

If $x_i = a + i \cdot d$, and we want $y'_{\sigma(0)},...,
y'_{\sigma(n)} $ to be an arithmetic progression for some $ \sigma
\in S_{n+1}$, then it must hold
$$
\begin{array}{lcr}
      x'_i = u^2 (a + id)+ s  &  & \\
      y'_i = b + \sigma(i) d',  &  & \mbox{for} \ \ i=0,..., n.\\
\end{array}
$$

We can take, with no loss of generality, $u=1$,
$s=0$, $t=0$. This involves only choosing an appropriate reference
system by translation and scaling (which would not affect s.a.p. in any case). 
Then we have
$$
\begin{array}{lcr}
      x'_i = x_i = a + id  &  & \\
      y'_i = y_i - r(a+id) = b + \sigma(i) d',   &  & \mbox{for} \ \ i=0,...,n.\\
\end{array}
$$
\noindent for some $b,d' \in K$.

These last identities can be written as a system of $n+1$ linear
equations in $r$, $b$ and $d'$, with matrix
$$
A^* = \left(
      \begin{array}{ccc|c}
        a+0d & 1 & \sigma(0) & y_0 \\
        a+1d & 1 & \sigma(1) & y_1 \\
        \vdots & \vdots & \vdots & \vdots \\
        a+nd & 1 & \sigma(n) & y_n \\
      \end{array}
\right)
$$

Note that the $y$--sequence $y_0,...,y_n$ is not an arithmetic
progression if and only if first, second and fourth columns are
independent; equivalently
$$
\exists s \in \{2,..., n\} \ \ \mbox{ such that} \ \ \begin{array}{|lcr|}
-0 & 1 & y_0 \\ -1 & 1 & y_1 \\ -s & 1 & y_s \\ \end{array} \neq 0.
$$

\vs\vs

\noindent {\bf Algorithm 1.--} Our first detection algorithm 
is based on the fact that, the existence of a solution to our
system (that is, the existence of a s.a.p.) is equivalent to $A^*$ 
having rank $3$. The formal algorithm goes like this:

\vs

\noindent {\bf Input Data:} $E,x_0,...,x_n$ (equivalently $E,x_0,n,d$).

\vs

\noindent {\bf Step 0:} Choose a suitable set $\{ y_0,...,y_n\}$ 
such that $(x_i,y_i) \in E(K)$.

\vs

\noindent {\bf Step 1:} (Fool--proof checking) Check whether 
$\{y_0,...,y_n\}$ is an arithmetic progression. If so, we are finished; 
if not, find an $s$, $2 \leq s \leq n$ as above. 

\vs

\noindent {\bf Step 2:} For any $\sigma \in S_n$ and any 
$i \in \{2,...,n\}$, $i \neq s$, compute the minor formed 
by the first, second, $s$--th and $i$--th rows of $A^*$.

\vs

\noindent {\bf Step 3:} If, for some $\sigma \in S_n$ the $n-2$ minors 
are null, solve the system to find $r$, $b$ and $d'$. If not, back to step 0.

\vs\vs

The main inconvenience of this procedure is its needing of $2^{n+1}(n+1)!(n-2)$ 
determinant computations, as there are $(n+1)!$ possibilities for $\sigma$ and
$2$ possibilities for each $y_i$. So we will try to find a more
efficient procedure, although this set--up will prove useful later on.

\vs\vs

\noindent {\bf Algorithm 2.--} Consider the affine points
$Q_0=(0, y_0, \sigma(0)),...,Q_n=(nd, y_n, \sigma(n)) \in {\bf A}^3(K)$. 
Note that, what we need from all them, in order to have a s.a.p., is to 
be in the same plane. This is the basis for our algorithm, whose input data and
steps 0 and 1 are identical to the previous one:

\vs

\noindent {\bf Step 2:} For each $\{i,j,k\} \subset \{0,1,...,
n\}$, we consider the plane
$$
\pi_{ijk} =  \langle (0,y_0,i),(d,y_1,j),(sd,y_s,k)\rangle,
$$

\vs

\noindent {\bf Step 3:} For $l=2,...,n$ and $l \neq s$, we
intersect the line
$$
x=ld, \quad  y=y_l
$$
\noindent with the plane $\pi_{ijk}$.

\vs

\noindent {\bf Step 4a:} If any of these intersections gives a
point $(ld, y_l, z_l)$ such that $z_l \notin \{0, 1,...,n\}$ or
$z_l$ is equal to another $z_{l'}$, then $\{z_0,...,z_n \}$ does
not correspond to $\{\sigma(0),...,\sigma(n) \}$ for any $\sigma
\in S_{n+1}$. Back to step 2, change the plane and repeat the
process or back to step 0 if all planes have been exhausted.

\vs

\noindent {\bf Step 4b:} If we find out a set of points $
Q_0=(x_0, y_0, z_0),..., Q_n=(x_n, y_n,  z_n)$ with $z_i =
\sigma(i)$ for $i=0,...,n$ and $\sigma \in S_{n+1}$, then $\sigma$ 
allows us to have a solution $r$, $b$, $d'$ of our system.

\vs\vs

As for computations is concerned, note that we have
$(n+1)n(n-1)/6$ possibilities for $\pi_{ijk}$, and for each plane
we have, at most, $n-2$ intersections. This, together with the
$2^{n+1}$ possibilities for each $y_i$, means a saving of around
$(n-3)!$ computations.

The implementation of both algorithms shows the time difference is not huge
(as $n$ dos not go very far), but it already grows significantly for $n\leq5$.

\vs

From now on, we will note $\sigma \in S_{n+1}$ by $\sigma = \left( a_0...a_n \right)$,
meaning $\sigma(0)=a_0,...,\sigma(n)=a_n$.

\vs

\noindent {\bf Example.--} Let the curve $Y^2 = X^3 -112 X + 400$, 
defined over $\QQ$, which has the following $x$--arithmetic progression of 
length $4$:
$$
x_k: -4, 0, 4, 8,
$$
\noindent and there are 4 of the $2^4$ $y$--sequences that lead to
simultaneous arithmetic progressions.
$$
\begin{array}{|rrrr|c|rrrr|}
        \hline
       y_0 & y_1 & y_2 & y_3 & \sigma & y'_0 & y'_1 & y'_2 & y'_3 \\
        \hline
        \hline
    28  & -20 &  4  &  4  &  (1023) & 4 & -20 & 28 & 52 \\ \hline
        &&&                 &  (0213) & -44/3 & -20 & -52/3 & -68/3 \\
    -28  & -20 &  -4  &  4  &  (1032) & -18 & -20 & -14 & -16 \\
        &&&                 &  (1302) & -84/5 & -20 & -76/5 & -92/5 \\ \hline
        &&&               &  (0213) & 44/3 & 20 & 52/3 & 68/3\\
    28  & 20 &  4  &  -4  &  (1032) &  18 & 20 & 14 & 16 \\
        &&&               &  (1302) & 84/5 & 20 & 76/5 & 92/5 \\ \hline
    -28  & 20 &  -4  &  -4  &  (1023) & -4 & 20 & -28 & -52 \\
       \hline
\end{array}
$$

Now this $x$--arithmetic progression can be extended to one
of length $5$:
$$
x_k: -4, 0, 4, 8, 12,
$$
\noindent and there are 2 of the $2^5$ possible $y$--sequences
that lead to simultaneous arithmetic progressions.

$$
\begin{array}{|rrrrr|c|rrrrr|}
        \hline
       y_0 & y_1 & y_2 & y_3 &y_4 & \sigma & y'_0 & y'_1 & y'_2 & y'_3 & y'_4 \\
        \hline
        \hline
    -28  & -20 &  -4  &  4  &  28 & (13240) & -44/3 & -20 & -52/3 & -68/3& -12 \\ \hline
    28  & 20 &  4  &  -4  &  -28 & (13240) & 44/3 & 20 & 52/3 & 68/3 & 12 \\
       \hline
\end{array}
$$

The equation for both cases is
$$
Y^2 - \frac{20}{3} XY = X^3 - \frac{100}{9}X^2-112X+400.
$$

Now, if we try to repeat the procedure for length $6$ with
$$
x_k: -4, 0, 4, 8, 12, 16,
$$
\noindent we find that none of the $2^6$ possible $y$--sequences
leads to a simultaneous arithmetic progression.

\vs\vs\vs

\noindent {\bf Open problem 1:} Find a procedure for deciding
whether an elliptic curve has an $x$--a.p. of given length.

\vs\vs

The most interesting results on this line are the parametrizations
by Bremner in \cite{BEM} which will be used afterwards in this
paper. However they are still far from being useful from a
computational point of view.

\section{The subsequence problem}

The programs developed in the previous section (specially the second one) 
were of great help with testing the examples we had created with the 
techniques shown on \cite{GT} and also with creating new ones. The 
counterexamples announced in Theorem 2 were product of these extensive
calculations. Here we present the simplest one.

\vs\vs

\noindent {\bf (Counter)Example.--} Consider the following
elliptic curve over $\QQ$, in Tate normal form,
$$
E \left( \frac{25}{21}, \frac{-2}{7} \right) : Y^2 + \frac{25}{21}
XY - \frac{2}{7} Y = X^3 +\frac{2}{7} X^2,
$$
\noindent which has the $x$--arithmetic progression of length $5$:
$$
x_k: \frac{-6}{7}, \frac{-4}{7}, \frac{-2}{7}, 0, \frac{2}{7}.
$$

Using the above procedure we find a $y$--sequence that gives
simultaneous arithmetic progression:
$$
y_k: \frac{4}{7}, \frac{16}{147}, \frac{92}{147}, 0, \frac{4}{21}.
$$

There is only one permutation $\sigma$ which passes Algorithm 2 and,
henceforth, allows the change of variables, $\sigma = (20413)$. The 
$y'$--arithmetic progression is
$$
y'_k: \frac{8}{49}, \frac{-8}{49}, \frac{24}{49}, 0, \frac{16}{49}.
$$
\noindent for the equation 
$$
Y^2+\frac{5}{21}XY - \frac{2}{7}Y = X^3+ \frac{92}{147}X^2-\frac{20}{147}X.
$$

This way, we have found a simultaneous arithmetic progression of
length $5$ not containing a simultaneous arithmetic progression of
length $4$, because the permutation involved is not the extension
of an $S_4$ permutation. In our many calculations these are 
singular cases: first of all, the permutation found is seldom
unique and, among the collected ones, we usually find an
extension of some $S_4$ permutation. But the fact is that these
counterexamples happen, which, by the way, carry an additional
difficulty for all arguments involving induction.

Interestingly enough, there are other suitable s.a.p. of length
$5$ with the same support (that is, other choices for the $y_i$) 
which happen to have subsequences of length $4$.

\vs\vs

\noindent {\bf Open problem 2:} Given an elliptic curve with a
s.a.p. of length $n$, is there always a s.a.p. of length $n-1$
whose support is contained in the support of the given one?

\vs\vs

We have found no examples to support a negative answer to this
question which, by the way, may serve as a weak induction result.

\section{The bound problem}

Our final look will be to the bound problem. As it was pointed out
by Bremner in \cite{BEM}, this kind of problems tend to become
unmanageable quite quickly. From our many attempts, we will
describe here the most successful of them all; which relies on a parametrization 
of curves with $x$--a.p. due to Bremner (\cite{BEM}) (here slightly changed 
for our purposes). In what follows we will assume $K=\QQ$. Note that all the
previous arguments do not rely on the base field at all.

\vs

First of all, we will parametrize elliptic curves in short Weierstrass form
$$
Y^2 = X^3+AX+B
$$
\noindent with four points in $x$--a.p.;
$$
P_0=(a,y_0), \; \; P_1=(a+d,y_1), \; \; P_2=(a+2d,y_2), \; \; P_3=(a+3d,y_3).
$$

Now, we consider the four polynomials $F_0,...,F_3$ 
given by
$$
F_i = y_i^2 - (a+id)^3 - A(a+id) - B,
$$
\noindent in $\QQ [y_0,y_1,y_2,y_3,a,d,A,B]$, and compute a Gr\"obner basis of the ideal
$\langle F_0,...,F_3 \rangle$. The \verb|tdeg| ordering in Maple V gives a 
basis of ten elements which can be used for computing $a,A,B$ taking $d,y_0,...,y_3$ 
as parameters:
\begin{eqnarray*}
A &=& \displaystyle \frac{-1}{6^2d^4} \left( y_0^4-9y_0^2y_1^2+6y_2^2y_0^2+y_3^2y_0^2+
21y_1^4-39y_2^2y_1^2+6y_3^2y_1^2+ \right. \\
&& \quad\quad \left. 21y_2^4-9y_3^2y_2^2+y_3^4 \right) = -P/36d^4 \\
B &=& \displaystyle \frac{1}{6^3d^6} \left( y_3^4y_0^2+4y_3^4y_1^2+y_3^4y_2^2-9y_3^2y_2^4 
-8y_3^2y_2^2y_0^2+24y_3^2y_1^4 \right. \\ 
&& \quad\quad -8y_3^2y_0^2y_1^2-12y_3^2y_2^2y_1^2+y_3^2y_0^4+y_0^4y_1^2-9y_0^2y_1^4+20y_1^6 \\ 
&& \quad\quad -21y_1^4y_2^2 +4y_0^4y_2^2+20y_2^6-21y_1^2y_2^4 +24y_0^2y_2^4 \\
&& \quad\quad \left. -12y_0^2y_2^2y_1^2 \right) = Q/6^3d^6 \\
a &=& \displaystyle \frac{-1}{6d^2} \left( -2y_0^2+5y_1^2-4y_2^2+y_3^2 \right) = -R/6d^2
\end{eqnarray*}
\noindent and, in addition, the first member of the basis is
$$
-y_3^2+6d^3+y_0^2-3y_1^2+3y_2^2.
$$

\vs

Now, making the scaling with $u=6d$ we obtain the following parametrization:
$$
\begin{array}{rl}
\mbox{Curve: } & Y^2 = X^3 - 6^2P + 6^3Q \\
\mbox{First term: } & 36d^2a = -6R \\
\mbox{Difference: } & 36d^3 = 6y_3^2-6y_0^2+18y_1^2-18y_2^2
\end{array}
$$

We will use from now on $A,B,a,d$ for these new polynomials. It is interesting noting
that the points in $x$--a.p. are now 
$$
P_i=(a+id, \,\pm 6y_id), \mbox{ for } i=0,...,3.
$$ 

\vs

We will try to produce curves with a s.a.p. of given length with a variant of 
Algorithm 1 which we will illustrate with the case of length $6$. 
In fact, using this procedure we might compute all curves with such a s.a.p. in 
contrast with the lack of exhaustiveness of \cite{GT}. If we want points 
$P_4=(a+4d,\, z_4)$ and $P_5=(a+5d,\, z_5)$ to be in the curve it must hold
\begin{eqnarray*}
z_4 &=& \pm 36\sqrt{4y_3^2+4y_1^2-y_0^2-6y_2^2}(-y_3^2+y_0^2-3y_1^2+3y_2^2) \\
z_5 &=& \pm 36\sqrt{-4y_0^2-20y_2^2+15y_1^2+10y_3^2}(-y_3^2+y_0^2-3y_1^2+3y_2^2)
\end{eqnarray*}

Hence, for the sake of consistency, we will call
\begin{eqnarray*}
y_4^2 &=& 4y_3^2+4y_1^2-y_0^2-6y_2^2 \\
y_5^2 &=& -4y_0^2-20y_2^2+15y_1^2+10y_3^2
\end{eqnarray*}
\noindent and our new points will then be $P_i=(a+id,\, \pm 6y_id)$ for $i=4,5$.

\vs

\noindent {\bf Example.--} We will show how to proceed using the $3$--cycle 
$\sigma =  (210345)$. After the above remarks, we may use the matrix $M$, given by
$$
M = \left( \begin{array}{ccc|c}
        0 & 1 & 1 & y_0 \\ 1 & 1 & 2 & y_1 \\ 2 & 1 & 0 & y_2 \\
        3 & 1 & 3 & y_3 \\ 4 & 1 & 4 & y_4 \\ 5 & 1 & 5 & y_5 \\
      \end{array}
\right),
$$
\noindent instead of the original matrix $A^*$ from Algorithm 1; and ask $M$ to 
have rank $3$, as we are assuming $d \neq 0$. Note that $\{y_0,...,y_5\}$ are 
{\em not} the $y$--coordinates of the points $P_0,...,P_5$. As all the minors are
linear polynomials on $y_0,...,y_5$ and we have also two quadratic relations, it is
not surprising that the complete solutions are two linear varieties, actually a 
plane and a line, given by the following parametrizations:
$$
\left\{y_0, \, \frac{y_4 + 3y_0}{4}, \, \frac{y_4+y_0}{2}, \, \frac{3y_4+y_0}{4}, \,
y_4, \, \frac{5y_4-y_0}{4} \right\},
$$
$$
\left\{y_0, \, \frac{7y_0}{19}, \, \frac{y_0}{19}, \, \frac{-15y_0}{19}, \, \frac{-27y_0}{19}, 
\, \frac{-39y_0}{19} \right\},
$$
\noindent from which the first one only contains points inducing $d=0$ and 
therefore must be discarded. In fact, these trivial solutions appear in all 
cases, which is clearly a by--product of our previous assumptions.

Now we make the substitutions induced by the second parametrization, obtaining
$$
A = \frac{-7840512}{130321}y_0^4, \quad B= \frac{8449090560}{47045881}y_0^6, \quad 
a = \frac{1536}{361}y_0^2, \quad d = \frac{48}{361}y_0^2,
$$
\noindent and the linear system given by $A^*$ has solution
$$
r= \frac{-60y_0}{19}, \quad b=\frac{98208y_0^3}{19^3}, \quad d'=\frac{-576y_0^3}{19^3}.
$$

This gives, after the corresponding substitution, the equation
$$
Y^2- \frac{120y_0}{19}XY = X^3 - \frac{3600y_0^2}{361}X^2 - \frac{7840512y_0^4}{130321}X 
+ \frac{8449090560y_0^6}{47045881},
$$
\noindent which has the following s.a.p.
$$
\left( \frac{1536y_0^2}{19^2}, \frac{97632y_0^3}{19^3} \right), \,
\left( \frac{1584y_0^2}{19^2}, \frac{97056y_0^3}{19^3} \right), \,
\left( \frac{1632y_0^2}{19^2}, \frac{978208y_0^3}{19^3} \right),
$$
$$
\left( \frac{1680y_0^2}{19^2}, \frac{96480y_0^3}{19^3} \right), \,
\left( \frac{1728y_0^2}{19^2}, \frac{95904y_0^3}{19^3} \right), \,
\left( \frac{1776y_0^2}{19^2}, \frac{95328y_0^3}{19^3} \right); \,
$$

All these curves are isomorphic to (an {\em easy}) one given by the case $y_0=19/2$;
$$
Y^2 - 60 YX = X^3 - 900 X^2 - 490032 X + 132017040, 
$$
having the sequence $\{ (384,12204)$, $(396,12132)$, $(408,12276)$,
$(420,12060)$, $(432,11988)$, $(444,11916)\}$.

\vs

We have not computed all curves with s.a.p. of length $6$, although we have 
bounded the number of curves by $19200$ cases, using the previous computations 
with all possible sign and permutation choices, counting only the number of 
possible solutions, that is, cases where the line does not induce $d=0$. 

\vs

To be precise, only half of the sign choices have to be 
considered, as every arithmetic progression of difference $d$ is also an 
arithemtic progression of difference $-d$, and hence, every curve appears 
at least twice, for a pair of inverse choice of signs and permutations.

Even so, not all of these cases are non--isomorphic elliptic curves; there 
might be isomorphic curves among them as well as genus $0$ curves. In the 
appendix we have given some explicit data for the first $100$ curves 
actually computed with this method where repeated curves already appear (in fact,
there are only $56$ non--isomorphic curves).

\vs

As a side remark, the distribution of the possible is extremely regular: 
there are $600$ allowed permutations (that is $600$ lines not giving $d=0$) 
for every sign choice and only sign changes were allowed most for the 
solutions for a fixed permutation. By the way, these differences ususally 
disappeared when finding the solutions to the system given by $A^*$. This, 
together with the repeated cases shown in the appendix, gives a heuristic 
estimation of only around $350$ non--isomorphic curves, but filling the 
details of such a list is beyond our computational possibilities so far.  

\vs

As for length $7$ is concerned our procedure shows there are no solutions, 
as all induce $d=0$. This case exhausts the possibilities of computer checking, 
at least with these methods, as it took around $20$ hours of CPU (which implied four days 
in real time) and, more constraining, $2^{11}$ Mb of stack memory. For considering 
this attack to length $8$, these figures should be multiplied at least by $16$ ($8$ for
the number of permutations and there are twice as many minors now), let alone the additional
difficulty of adding a new quadratic polynomial to the system, which is not easy to
measure.

\vs\vs\vs

\noindent {\bf Open problem 3:} Find a universal bound for the length of a s.a.p.
on elliptic curves over $\QQ$.

\vs\vs

Note that an affirmative answer to the open problem $2$ would mean $6$ is
the answer to open problem $3$.

\section*{Appendix: Examples of curves with s.a.p. of length $6$}

Please note that this stream was computer--generated taking as a unique choice
the permutation. Due to this, some of the curves (for instance examples $009$ 
and $010$) appear more than once, as they have different s.a.p. Also one may 
find isomorphic curves (as $001$ and $002$). As we said above, only $56$ 
non--isomorphic curves can be found in this table, but we have preferred to
leave as it came, as we feel it illustrated better the phenomenon. The entries 
of the table after each equation are:

\vs

1) Permutation ($\sigma$): noted as above by $(\sigma(0)\sigma(1)...\sigma(5))$.

\vs

2) Numerical data (N.D.): The set $[a,d,b,d']$ which fits the equation.

\vs

3) Rank ($r$): The rank of the curve, computed with APECS ($\leq$ 
means APECS failed to actually compute the rank, in which case the best upper 
bound given is shown). All curves have trivial torsion group.

{\scriptsize\begin{eqnarray*}
\hline &&\\
001 &&Y^2- 180 YX + 8100 X^2 - X^3 + 4892251392 X- 134063884477440 \\
        && \sigma = (321450),\; \mbox{N.D.}=[66432, -20304, 13044672, -5725728], \; r=5 \\
&& \\ \hline  &&\\
002 &&Y^2+ 180 YX + 8100 X^2 - X^3 + 4892251392 X- 134063884477440 \\
        && \sigma = (054123),\; \mbox{N.D.}=[-35088, 20304, -13044672, 5725728], \; r=5 \\
&& \\ \hline  &&\\
003 &&Y^2+ 20 YX + 100 X^2 - X^3 + 36478512 X- 82321246080 \\
        && \sigma = (423150),\; \mbox{N.D.}=[5724, -1584, -367704, 139392], \; r=5 \\
&& \\ \hline  &&\\
004 &&Y^2- 60 YX + 900 X^2 - X^3 + 86832 X- 8864640 \\
        && \sigma = (534201),\; \mbox{N.D.}=[-324, 144, 9288, -3456], \; r=2 \\
&& \\ \hline  &&\\
005 &&Y^2- 20 YX + 100 X^2 - X^3 + 466992 X- 549797760 \\
        && \sigma = (045312),\; \mbox{N.D.}=[1308, -432, -33576, 3456], \; r=3 \\
&& \\ \hline  &&\\
006 &&Y^2- 2 YX + X^2 - X^3 + 238707 X- 41709006 \\
        && \sigma = (150423),\; \mbox{N.D.}=[-513, 180, 8487, -3600], \; r=3 \\
&& \\ \hline  &&\\
007 &&Y^2- 180 YX + 8100 X^2 - X^3 + 51432192 X- 368371860480 \\
        && \sigma = (034521),\; \mbox{N.D.}=[-9312, 4752, -638496, 313632], \; r=4 \\
&& \\ \hline  &&\\
008 &&Y^2- 612 YX + 93636 X^2 - X^3 + 26962612992 X- 1882111797863424 \\
        && \sigma = (145032),\; \mbox{N.D.}=[-164832, 66960, -19904832, 12454560], \; r=5 \\
&& \\ \hline  &&\\
009 &&Y^2+ 396 YX + 39204 X^2 - X^3 + 388106595072 X- 88686989876929536 \\
        && \sigma = (250143),\; \mbox{N.D.}=[-715488, 235440, 459510624, -153977760], \; r=5 \\
&& \\ \hline  &&\\
010 &&Y^2- 576 YX + 82944 X^2 - X^3 + 271059091200 X- 49072046238950400 \\
        && \sigma = (134502),\; \mbox{N.D.}=[734160, -265680, -415035360, 130714560], \; r=5 \\
&& \\ \hline  &&\\
011 &&Y^2- 576 YX + 82944 X^2 - X^3 + 271059091200 X- 49072046238950400 \\
        && \sigma = (350124),\; \mbox{N.D.}=[-594240, 265680, 238537440, -130714560], \; r=5 \\
&& \\ \hline  &&\\
012 &&Y^2+ 396 YX + 39204 X^2 - X^3 + 388106595072 X- 88686989876929536 \\
        && \sigma = (305412),\; \mbox{N.D.}=[-715488, 235440, -310378176, 153977760], \; r=5 \\
&& \\ \hline  &&\\
013 &&Y^2- 612 YX + 93636 X^2 - X^3 + 26962612992 X- 1882111797863424 \\
        && \sigma = (410523),\; \mbox{N.D.}=[-164832, 66960, 42367968, -12454560], \; r=5 \\
&& \\ \hline  &&\\
014 &&Y^2- 180 YX + 8100 X^2 - X^3 + 51432192 X- 368371860480 \\
        && \sigma = (521034),\; \mbox{N.D.}=[-9312, 4752, 929664, -313632], \; r=4 \\
&& \\ \hline  &&\\
015 &&Y^2- 300 YX + 22500 X^2 - X^3 + 894074112 X- 19561912750080 \\
        && \sigma = (230145),\; \mbox{N.D.}=[-30384, 11808, -6045408, 2904768], \; r=6 \\
&& \\ \hline  &&\\
016 &&Y^2- 60 YX + 900 X^2 - X^3 + 48045312 X- 189087436800 \\
        && \sigma = (341250),\; \mbox{N.D.}=[-7344, 3744, -714528, 292032], \; r=5 \\
&& \\ \hline  &&\\
017 &&Y^2- 60 YX + 900 X^2 - X^3 + 48045312 X- 189087436800 \\
        && \sigma = (503412),\; \mbox{N.D.}=[11376, -3744, 745632, -292032], \; r=5 \\
&& \\ \hline  &&\\
018 &&Y^2+ 300 YX + 22500 X^2 - X^3 + 894074112 X- 19561912750080 \\
        && \sigma = (014523),\; \mbox{N.D.}=[28656, -11808, -8478432, 2904768], \; r=6 \\
&& \\ \hline  &&\\
019 &&Y^2- 1800 YX + 810000 X^2 - X^3 + 7410269014272 X- 8108956611489899520 \\
        && \sigma = (302145),\; \mbox{N.D.}=[2763696, -577584, 3502719072, -1323822528], \; r\leq 6 \\
&& \\ \hline  &&\\
020 &&Y^2+ 2880 YX + 2073600 X^2 - X^3 + 403246536682752 X- 3207719147336296058880 \\
        && \sigma = (413250),\; \mbox{N.D.}=[-829968, 5312736, 62302243392, -30155089536], \; r\leq 6 \\
&& \\ \hline  &&\\
021 &&Y^2+ 4680 YX + 5475600 X^2 - X^3 + 3994541989632 X- 6207305561351930880 \\
        && \sigma = (524301),\; \mbox{N.D.}=[-2327712, 728784, -3390388704, 2107643328], \; r\leq 7 \\
&& \\ \hline  &&\\
022 &&Y^2- 2880 YX + 2073600 X^2 - X^3 + 8792422324992 X- 14409886139859502080 \\
        && \sigma = (035412),\; \mbox{N.D.}=[7563792, -2228688, -8618074272, 1791865152], \; r\leq 11 \\
&& \\ \hline  &&\\
023 &&Y^2+ 7632 YX + 14561856 X^2 - X^3 + 1208964238739712 X- 23944945010158503235584 \\
        && \sigma = (140523),\; \mbox{N.D.}=[-42040752, 12290400, 265952884032, -83918851200], \; r\leq 7 \\
&& \\ \hline  &&\\
024 &&Y^2- 6912 YX + 11943936 X^2 - X^3 + 214026269587200 X- 1973672371225801958400 \\
        && \sigma = (250134),\; \mbox{N.D.}=[-9752640, 7294320, 76016037600, -26872274880], \; r\leq 8 \\
&& \\ \hline  &&\\
025 &&Y^2+ 120 YX + 3600 X^2 - X^3 + 81948672 X- 327577374720 \\
        && \sigma = (532140),\; \mbox{N.D.}=[8064, -3168, 1169856, -418176], \; r=4 \\
&& \\ \hline  &&\\
026 &&Y^2+ 960 YX + 230400 X^2 - X^3 + 92277352334592 X- 342160070800370356224 \\
        && \sigma = (043251),\; \mbox{N.D.}=[11544048, -3013920, -34095021120, 10813944960], \; r=4 \\
&& \\ \hline  &&\\
027 &&Y^2+ 2616 YX + 1710864 X^2 - X^3 + 3327107104512 X- 3360005144504534016 \\
        && \sigma = (154302),\; \mbox{N.D.}=[-515808, 403920, -2454087456, 906396480], \; r\leq 6 \\
&& \\ \hline  &&\\
028 &&Y^2- 3072 YX + 2359296 X^2 - X^3 + 2328753038487552 X- 49786655684722942869504 \\
        && \sigma = (205413),\; \mbox{N.D.}=[-55130112, 19313280, 240362650368, -111012733440], \; r\leq 6 \\
&& \\ \hline  &&\\
029 &&Y^2- 4944 YX + 6110784 X^2 - X^3 + 14552162687232 X- 21391150302252205056 \\
        && \sigma = (310524),\; \mbox{N.D.}=[1164672, 340560, 5468006304, -44953920], \; r\leq 7 \\
&& \\ \hline  &&\\
030 &&Y^2- 960 YX + 230400 X^2 - X^3 + 19321742592 X- 1033886822584320 \\
        && \sigma = (421035),\; \mbox{N.D.}=[60912, 6048, 38636352, -72576], \; r=5 \\
&& \\ \hline  &&\\
031 &&Y^2+ 40 YX + 400 X^2 - X^3 + 8854272 X- 8259978240 \\
        && \sigma = (153420),\; \mbox{N.D.}=[-2448, 1584, 242208, -69696], \; r=4 \\
&& \\ \hline  &&\\
032 &&Y^2- 64 YX + 1024 X^2 - X^3 + 1389312 X- 394896384 \\
        && \sigma = (204531),\; \mbox{N.D.}=[-1152, 720, -44064, 14400], \; r=3 \\
&& \\ \hline  &&\\
033 &&Y^2- 64 YX + 1024 X^2 - X^3 + 1389312 X- 394896384 \\
        && \sigma = (420153),\; \mbox{N.D.}=[2448, -720, 27936, -14400], \; r=3 \\
&& \\ \hline  &&\\
034 &&Y^2- 40 YX + 400 X^2 - X^3 + 8854272 X- 8259978240 \\
        && \sigma = (531204),\; \mbox{N.D.}=[5472, -1584, 106272, -69696], \; r=4 \\
&& \\ \hline  &&\\
035 &&Y^2- 1080 YX + 291600 X^2 - X^3 + 4379037151491072 X- 136063103467431710822400 \\
        && \sigma = (240315),\; \mbox{N.D.}=[-77988192, 33054048, -438167226816, 226486336896], \; r\leq 8 \\
&& \\ \hline  &&\\
036 &&Y^2- 720 YX + 129600 X^2 - X^3 + 9969629720832 X- 13778174775900128256 \\
        && \sigma = (351420),\; \mbox{N.D.}=[-3640128, 1648080, -6422790240, 2155688640], \; r=1 \\
&& \\ \hline  &&\\
037 &&Y^2+ 5544 YX + 7683984 X^2 - X^3 + 422337203926272 X- 3401406251614088487936 \\
        && \sigma = (402531),\; \mbox{N.D.}=[-18956688, 10568880, 102961496736, -29550588480], \; r\leq 7 \\
&& \\ \hline  &&\\
038 &&Y^2+ 1080 YX + 291600 X^2 - X^3 + 4379037151491072 X- 136063103467431710822400 \\
        && \sigma = (513042),\; \mbox{N.D.}=[87282048, -33054048, 438167226816, -226486336896], \; r\leq 8 \\
&& \\ \hline  &&\\
039 &&Y^2- 720 YX + 129600 X^2 - X^3 + 9969629720832 X- 13778174775900128256 \\
        && \sigma = (024153),\; \mbox{N.D.}=[4600272, -1648080, -6422790240, 2155688640], \; r=1 \\
&& \\ \hline  &&\\
040 &&Y^2- 5544 YX + 7683984 X^2 - X^3 + 422337203926272 X- 3401406251614088487936 \\
        && \sigma = (135204),\; \mbox{N.D.}=[33887712, -10568880, -102961496736, 29550588480], \; r\leq 7 \\
&& \\ \hline  &&\\
041 &&Y^2+ 20 YX + 100 X^2 - X^3 + 466992 X- 549797760 \\
        && \sigma = (342015),\; \mbox{N.D.}=[-852, 432, 16296, 3456], \; r=3 \\
&& \\ \hline  &&\\
042 &&Y^2- 60 YX + 900 X^2 - X^3 + 86832 X- 8864640 \\
        && \sigma = (453120),\; \mbox{N.D.}=[396, -144, -7992, 3456], \; r=2 \\
&& \\ \hline  &&\\
043 &&Y^2- 20 YX + 100 X^2 - X^3 + 36478512 X- 82321246080 \\
        && \sigma = (504231),\; \mbox{N.D.}=[-2196, 1584, -329256, 139392], \; r=5 \\
&& \\ \hline  &&\\
044 &&Y^2- 4 YX + 4 X^2 - X^3 + 3819312 X- 2669376384 \\
        && \sigma = (231504),\; \mbox{N.D.}=[1548, -720, -76104, 28800], \; r=3 \\
&& \\ \hline  &&\\
045 &&Y^2- 2880 YX + 2073600 X^2 - X^3 + 110993807215872 X- 445121432212494274560 \\
        && \sigma = (402315),\; \mbox{N.D.}=[14071152, -3627936, 35652163392, -14061879936], \; r\leq 7 \\
&& \\ \hline  &&\\
046 &&Y^2+ 360 YX + 32400 X^2 - X^3 + 27909792000 X- 1526822144640000 \\
        && \sigma = (513420),\; \mbox{N.D.}=[-99600, 75600, 73720800, -22680000], \; r\leq 5 \\
&& \\ \hline  &&\\
047 &&Y^2- 7056 YX + 12446784 X^2 - X^3 + 129437671756032 X- 543091962160151110656 \\
        && \sigma = (135042),\; \mbox{N.D.}=[17527152, -4944240, -5119105824, 6467065920], \; r\leq 5 \\
&& \\ \hline  &&\\
048 &&Y^2- 7488 YX + 14017536 X^2 - X^3 + 841545692773632 X- 15278101189303210758144 \\
        && \sigma = (240153),\; \mbox{N.D.}=[-34464912, 10864800, -202374030528, 65579932800], \; r\leq 8 \\
&& \\ \hline  &&\\
049 &&Y^2- 6696 YX + 11209104 X^2 - X^3 + 1767944739520512 X- 27067552170879621206016 \\
        && \sigma = (351204),\; \mbox{N.D.}=[-46912128, 22001760, 291906110016, -122241778560], \; r\leq 7 \\
&& \\ \hline  &&\\
050 &&Y^2- 90 YX + 2025 X^2 - X^3 + 3434834187 X- 77407635019590 \\
        && \sigma = (053421),\; \mbox{N.D.}=[39747, -2376, 3706047, -627264], \; r=5 \\
&& \\ \hline  &&\\
051 &&Y^2+ 306 YX + 23409 X^2 - X^3 + 338359707 X- 3186238861494 \\
        && \sigma = (104532),\; \mbox{N.D.}=[-21237, 11880, 3783861, -1425600], \; r=4 \\
&& \\ \hline  &&\\
052 &&Y^2- 612 YX + 93636 X^2 - X^3 + 5413755312 X- 203919287135616 \\
        && \sigma = (320154),\; \mbox{N.D.}=[152652, -47520, 26753112, -11404800], \; r=4 \\
&& \\ \hline  &&\\
053 &&Y^2- 90 YX + 2025 X^2 - X^3 + 3434834187 X- 77407635019590 \\
        && \sigma = (431205),\; \mbox{N.D.}=[27867, 2376, 569727, 627264], \; r=5 \\
&& \\ \hline  &&\\
054 &&Y^2+ 396 YX + 39204 X^2 - X^3 + 388106595072 X- 88686989876929536 \\
        && \sigma = (214503),\; \mbox{N.D.}=[461712, -235440, -310378176, 153977760], \; r=5 \\
&& \\ \hline  &&\\
055 &&Y^2- 612 YX + 93636 X^2 - X^3 + 26962612992 X- 1882111797863424 \\
        && \sigma = (325014),\; \mbox{N.D.}=[169968, -66960, 42367968, -12454560], \; r=5 \\
&& \\ \hline  &&\\
056 &&Y^2- 180 YX + 8100 X^2 - X^3 + 51432192 X- 368371860480 \\
        && \sigma = (430125),\; \mbox{N.D.}=[14448, -4752, 929664, -313632], \; r=4 \\
&& \\ \hline  &&\\
057 &&Y^2- 48 YX + 576 X^2 - X^3 + 9262512 X- 8383430016 \\
        && \sigma = (250341),\; \mbox{N.D.}=[-3348, 1440, -123552, 43200], \; r=2 \\
&& \\ \hline  &&\\
058 &&Y^2- 24 YX + 144 X^2 - X^3 + 578907 X- 130991094 \\
        && \sigma = (412503),\; \mbox{N.D.}=[963, -360, 11556, -5400], \; r=2 \\
&& \\ \hline  &&\\
059 &&Y^2- 60 YX + 900 X^2 - X^3 + 767232 X- 253808640 \\
        && \sigma = (523014),\; \mbox{N.D.}=[1008, -288, 16416, -1728], \; r=3 \\
&& \\ \hline  &&\\
060 &&Y^2- 60 YX + 900 X^2 - X^3 + 767232 X- 253808640 \\
        && \sigma = (145230),\; \mbox{N.D.}=[-432, 288, 7776, 1728], \; r=3 \\
&& \\ \hline  &&\\
061 &&Y^2- 5184 YX + 6718464 X^2 - X^3 + 142783615579392 X- 708085096419827994624 \\
        && \sigma = (352041),\; \mbox{N.D.}=[13189008, -5387040, -40963299264, 13898563200], \; r=4 \\
&& \\ \hline  &&\\
062 &&Y^2- 2142 YX + 1147041 X^2 - X^3 + 15808640228307 X- 39293749322977053294 \\
        && \sigma = (403152),\; \mbox{N.D.}=[5885823, -1961820, 10341141993, -4072738320], \; r\leq 7 \\
&& \\ \hline  &&\\
063 &&Y^2- 612 YX + 93636 X^2 - X^3 + 52888734000 X- 42522157354464000 \\
	&& \sigma = (514203),\; \mbox{N.D.}=[-243780, 110160, 98574840, -44945280], \; r\leq 5 \\
&& \\ \hline  &&\\
064 &&Y^2- 396 YX + 39204 X^2 - X^3 + 16655166000 X- 690091901769600 \\
        && \sigma = (025314),\; \mbox{N.D.}=[250140, -71280, -60813720, 18817920], \; r=4 \\
&& \\ \hline  &&\\
065 &&Y^2- 1710 YX + 731025 X^2 - X^3 + 1041760268307 X- 785868626594613294 \\
        && \sigma = (130425),\; \mbox{N.D.}=[-1270977, 427140, -1427543055, 579201840], \; r\leq 5 \\
&& \\ \hline  &&\\
066 &&Y^2- 3240 YX + 2624400 X^2 - X^3 + 2359518446592 X- 1392039372309073920 \\
        && \sigma = (241530),\; \mbox{N.D.}=[502176, 175392, 1361800512, 14732928], \; r=7 \\
&& \\ \hline  &&\\
067 &&Y^2- 60 YX + 900 X^2 - X^3 + 48045312 X- 189087436800 \\
        && \sigma = (503412),\; \mbox{N.D.}=[11376, -3744, 745632, -292032], \; r=5 \\
&& \\ \hline  &&\\
068 &&Y^2+ 300 YX + 22500 X^2 - X^3 + 894074112 X- 19561912750080 \\
        && \sigma = (014523),\; \mbox{N.D.}=[28656, -11808, -8478432, 2904768], \; r=6 \\
&& \\ \hline  &&\\
069 &&Y^2- 300 YX + 22500 X^2 - X^3 + 894074112 X- 19561912750080 \\
        && \sigma = (230145),\; \mbox{N.D.}=[-30384, 11808, -6045408, 2904768], \; r=6 \\
&& \\ \hline  &&\\
070 &&Y^2- 60 YX + 900 X^2 - X^3 + 48045312 X- 189087436800 \\
        && \sigma = (341250),\; \mbox{N.D.}=[-7344, 3744, -714528, 292032], \; r=5 \\
&& \\ \hline  &&\\
071 &&Y^2- 960 YX + 230400 X^2 - X^3 + 13201887036672 X- 17372103328571019264 \\
        && \sigma = (502341),\; \mbox{N.D.}=[5909328, -1703520, 8039027520, -3454738560], \; r=6 \\
&& \\ \hline  &&\\
072 &&Y^2+ 600 YX + 90000 X^2 - X^3 + 4185188352 X- 119899753943040 \\
        && \sigma = (013452),\; \mbox{N.D.}=[52704, -12384, -22572864, 6390144], \; r\leq 7 \\
&& \\ \hline  &&\\
073 &&Y^2+ 768 YX + 147456 X^2 - X^3 + 28465015272192 X- 52273748736104702976 \\
        && \sigma = (124503),\; \mbox{N.D.}=[4297392, -2075760, -8562479328, 3860913600], \; r\leq 6 \\
&& \\ \hline  &&\\
074 &&Y^2- 4368 YX + 4769856 X^2 - X^3 + 102164509274112 X- 521183246173070966784 \\
        && \sigma = (235014),\; \mbox{N.D.}=[7567488, -3850560, 20709101952, -8825483520], \; r\leq 6 \\
&& \\ \hline  &&\\
075 &&Y^2- 192 YX + 9216 X^2 - X^3 + 1067634432 X- 20407326188544 \\
        && \sigma = (340125),\; \mbox{N.D.}=[-32832, 12240, -6162912, 2496960], \; r=4 \\
&& \\ \hline  &&\\
076 &&Y^2+ 120 YX + 3600 X^2 - X^3 + 210573084672 X- 33546735732363264 \\
        && \sigma = (451230),\; \mbox{N.D.}=[-495072, 211680, -339888960, 124467840], \; r=6 \\
&& \\ \hline  &&\\
077 &&Y^2- 180 YX + 8100 X^2 - X^3 + 51432192 X- 368371860480 \\
        && \sigma = (125430),\; \mbox{N.D.}=[14448, -4752, -638496, 313632], \; r=4 \\
&& \\ \hline  &&\\
078 &&Y^2- 612 YX + 93636 X^2 - X^3 + 26962612992 X- 1882111797863424 \\
        && \sigma = (230541),\; \mbox{N.D.}=[169968, -66960, -19904832, 12454560], \; r=5 \\
&& \\ \hline  &&\\
079 &&Y^2+ 396 YX + 39204 X^2 - X^3 + 388106595072 X- 88686989876929536 \\
        && \sigma = (341052),\; \mbox{N.D.}=[461712, -235440, 459510624, -153977760], \; r=5 \\
&& \\ \hline  &&\\
080 &&Y^2- 60 YX + 900 X^2 - X^3 + 48045312 X- 189087436800 \\
        && \sigma = (214305),\; \mbox{N.D.}=[-7344, 3744, 745632, -292032], \; r=5 \\
&& \\ \hline  &&\\
081 &&Y^2- 300 YX + 22500 X^2 - X^3 + 894074112 X- 19561912750080 \\
        && \sigma = (325410),\; \mbox{N.D.}=[-30384, 11808, 8478432, -2904768], \; r=6 \\
&& \\ \hline  &&\\
082 &&Y^2+ 300 YX + 22500 X^2 - X^3 + 894074112 X- 19561912750080 \\
        && \sigma = (541032),\; \mbox{N.D.}=[28656, -11808, 6045408, -2904768], \; r=6 \\
&& \\ \hline  &&\\
083 &&Y^2- 60 YX + 900 X^2 - X^3 + 48045312 X- 189087436800 \\
        && \sigma = (052143),\; \mbox{N.D.}=[11376, -3744, -714528, 292032], \; r=5 \\
&& \\ \hline  &&\\
084 &&Y^2- 3240 YX + 2624400 X^2 - X^3 + 2359518446592 X- 1392039372309073920 \\
        && \sigma = (314025),\; \mbox{N.D.}=[502176, 175392, 1435465152, -14732928], \; r=7 \\
&& \\ \hline  &&\\
085 &&Y^2- 1710 YX + 731025 X^2 - X^3 + 1041760268307 X- 785868626594613294 \\
        && \sigma = (425130),\; \mbox{N.D.}=[-1270977, 427140, 1468466145, -579201840], \; r\leq 5 \\
&& \\ \hline  &&\\
086 &&Y^2- 396 YX + 39204 X^2 - X^3 + 16655166000 X- 690091901769600 \\
        && \sigma = (530241),\; \mbox{N.D.}=[250140, -71280, 33275880, -18817920], \; r=4 \\
&& \\ \hline  &&\\
087 &&Y^2- 612 YX + 93636 X^2 - X^3 + 52888734000 X- 4252215735446400 \\
        && \sigma = (041352),\; \mbox{N.D.}=[-243780, 110160, -126151560, 44945280], \; r\leq 5 \\
&& \\ \hline  &&\\
088 &&Y^2+ 4284 YX + 4588164 X^2 - X^3 + 252938243652912 X- 2514799956670531410816 \\
        && \sigma = (152403),\; \mbox{N.D.}=[23543292, -7847280, 80180396856, -32581906560], \; r\leq 7 \\
&& \\ \hline  &&\\
089 &&Y^2- 5184 YX + 6718464 X^2 - X^3 + 142783615579392 X- 708085096419827994624 \\
        && \sigma = (203514),\; \mbox{N.D.}=[13189008, -5387040, 28529516736, -13898563200], \; r=4 \\
&& \\ \hline  &&\\
090 &&Y^2- 60 YX + 900 X^2 - X^3 + 767232 X- 253808640 \\
        && \sigma = (410325),\; \mbox{N.D.}=[-432, 288, 16416, -1728], \; r=3 \\
&& \\ \hline  &&\\
091 &&Y^2- 60 YX + 900 X^2 - X^3 + 767232 X- 253808640 \\
        && \sigma = (032541),\; \mbox{N.D.}=[1008, -288, 7776, 1728], \; r=3 \\
&& \\ \hline  &&\\
092 &&Y^2- 24 YX + 144 X^2 - X^3 + 578907 X- 130991094 \\
        && \sigma = (143052),\; \mbox{N.D.}=[963, -360, -15444, 5400], \; r=2 \\
&& \\ \hline  &&\\
093 &&Y^2- 48 YX + 576 X^2 - X^3 + 9262512 X- 8383430016 \\
        && \sigma = (305214),\; \mbox{N.D.}=[-3348, 1440, 92448, -43200], \; r=2 \\
&& \\ \hline  &&\\
094 &&Y^2+ 120 YX + 3600 X^2 - X^3 + 81948672 X- 327577374720 \\
        && \sigma = (514320),\; \mbox{N.D.}=[-7776, 3168, -921024, 418176], \; r=4 \\
&& \\ \hline  &&\\
095 &&Y^2- 960 YX + 230400 X^2 - X^3 + 19321742592 X- 1033886822584320 \\
        && \sigma = (025431),\; \mbox{N.D.}=[91152, -6048, 38273472, 72576], \; r=5 \\
&& \\ \hline  &&\\
096 &&Y^2- 4944 YX + 6110784 X^2 - X^3 + 14552162687232 X- 21391150302252205056 \\
        && \sigma = (130542),\; \mbox{N.D.}=[2867472, -340560, 5243236704, 44953920], \; r\leq 7 \\
&& \\ \hline  &&\\
097 &&Y^2- 3072 YX + 2359296 X^2 - X^3 + 2328753038487552 X- 49786655684722942869504 \\
        && \sigma = (241053),\; \mbox{N.D.}=[41436288, -19313280, -314701016832, 111012733440], \; r\leq 6 \\
&& \\ \hline  &&\\
098 &&Y^2- 2616 YX + 1710864 X^2 - X^3 + 3327107104512 X- 3360005144504534016 \\
        && \sigma = (352104),\; \mbox{N.D.}=[1503792, -403920, -2077894944, 906396480], \; r\leq 6 \\
&& \\ \hline  &&\\
099 &&Y^2- 960 YX + 230400 X^2 - X^3 + 92277352334592 X- 342160070800370356224 \\
        && \sigma = (403215),\; \mbox{N.D.}=[-3525552, 3013920, -19974703680, 10813944960], \; r\leq 4 \\
&& \\ \hline  &&\\
100 &&Y^2- 60 YX + 900 X^2 - X^3 + 767232 X- 253808640 \\
        && \sigma = (145230),\; \mbox{N.D.}=[-432, 288, 7776, 1728], \; r=3 \\
&& \\ \hline  &&\\
\end{eqnarray*}}

As Bremner noticed in \cite{BEM}, points in arithmetic progression seem to have 
a tendency to be independent. If we only take into account the ranks actually 
computed, we get an average over $4$, when the average for random curves is known 
to be much smaller (\cite{Young}). 

It should be noticed that example $039$ is a remarkable case: no $x$--a.p. of 
length $6$ were known in groups of rank one (see the last remark of \cite{BEM}), 
and \cite{BST} shows why this is so uncommon. The relation, if any, between the 
rank and the length of s.a.p. (or $x$--a.p.) seems a much harder problem to tackle.

\vs\vs

Irene Garc\'{\i}a--Selfa (Email: \verb|igselfa@us.es|).

Jos\'e  M. Tornero (Corresponding author, email: \verb|tornero@us.es|).

\vs

Departamento de \'Algebra 

Facultad de Matem\'aticas

Universidad de Sevilla

Apdo. 1160. 41080 Sevilla (Spain).


\begin{thebibliography}{99}

\bibitem{BAA1}
A. Bremner: {\em On squares of squares.} Acta Arith. {\bf
LXXXVIII} (1999) 289--297.

\bibitem{BAA2}
A. Bremner: {\em On squares of squares II.} Acta Arith. {\bf XCIX}
(3) (2001) 289--308.

\bibitem{BEM}
A. Bremner: {\em On arithmetic progressions on elliptic curves.}
Experiment. Math. {\bf 8} (4) (1999) 409--413.
-
\bibitem{BST}
A. Bremner, J.H. Silverman, N. Tzanakis: {\em Integral points in
arithmetic progression on $y^2=x(x^2-n^2)$.} J. Number Theory {\bf
80} (2000) 187--208.

\bibitem{Campbell}
G. Campbell: {\em A note on arithmetic progressions on elliptic
curves.} J. Integer Seq. {\bf 6} (2003) 03.1.3.

\bibitem{GT}
I. Garc\'{\i}a--Selfa, J.M. Tornero: {\em Searching for
simultaneous arithmetic progressions on elliptic curves.} Bull.
Austral. Math. Soc. {\bf 71} (2005) 417--424.

\bibitem{LV}
J.B. Lee, W.Y. V\'elez: {\em Integral solutions in arithmetic
progression for $y^2=x^3+k$.} Period. Math. Hungar. {\bf 25}
(1992) 31--49.

\bibitem{Mohanty}
S.P. Mohanty: {\em On consecutive integral solutions for $y^2 =
x^3+k$.} Proc. Amer. Math. Soc. {\bf 48} (1975) 281--285.

\bibitem{Silv}
J.H. Silverman: {\em The arithmetic of elliptic curves.} Springer
(1986).

\bibitem{Young}
M.P. Young: {\em Low-lying zeros of families of elliptic curves.} 
 J. Amer. Math. Soc. {\bf 19} (2006) 205--250.
\end{thebibliography}
\end{document}